\documentclass[twoside,notitlepage,11pt]{article}

\usepackage{amssymb}
\usepackage[leqno]{amsmath}
\usepackage{amsfonts}
\usepackage{amsopn}
\usepackage{amstext}
\usepackage{amsthm}
\usepackage{hyperref}

\newcommand{\Aaa}{{\cal{A}}}
\newcommand{\Bee}{{\cal{B}}}
\newcommand{\Cee}{{\cal{C}}}

\newcommand{\Yu}{{\cal{U}}}

\newcommand{\Qyu}{{\Bbb{Q}}}
\newcommand{\Err}{{\Bbb{R}}}

\newcommand{\al}{\alpha}

\renewcommand{\phi}{\varphi}
\renewcommand{\rho}{\varrho}

\newcommand{\rest}{\restriction}

\newcommand{\loe}{\le}
\newcommand{\goe}{\ge}

\newcommand{\subs}{\subseteq}

\newcommand{\conv}{\operatorname{conv}}

\newcommand{\id}[1]{\operatorname{id}_{#1}}

\newcommand{\dom}{\operatorname{dom}}

\newtheorem{tw}{Theorem}[section]

\newtheorem{lm}[tw]{Lemma}
\newtheorem{prop}[tw]{Proposition}

\theoremstyle{definition}

\newtheorem{ex}[tw]{Example}
\theoremstyle{remark}

\newcommand{\setof}[2]{\{#1\colon #2\}}

\newcommand{\sett}[2]{\{#1\}_{#2}}
\newcommand{\sn}[1]{\{#1\}} 
\newcommand{\dn}[2]{\{#1,#2\}} 
\newcommand{\pair}[2]{\langle #1, #2 \rangle} 
\newcommand{\map}[3]{#1\colon #2 \to #3} 

\newcommand{\wek}[1]{{\vec{#1}}}

\newcommand{\tc}[1]{\operatorname{tc}\left(#1\right)}

\providecommand{\cal}{\mathcal}
\renewcommand{\Bbb}{\mathbb}

\newenvironment{pf}{\begin{proof}}{\end{proof}}

\newcommand{\suppt}{\operatorname{suppt}}

\newcommand{\xk}[1]{{\mathbb X}\left({#1}\right)}
\newcommand{\kx}[1]{{\mathbb K}\left({#1}\right)}

\newcommand{\cmp}{\circ} 

\newtheorem{question}{Question}

\title{The structure  of Valdivia compact lines}
\author{
{\sc Ond\v rej F.K. Kalenda}
\footnote{Research supported by Research project MSM 0021620839 and partly by Research grant GA\v{C}R 201/00/0018.}
\\
{\small Charles University, Faculty of Mathematics and Physics}\\
{\small Department of Mathematical Analysis}\\
{\small Sokolovsk\'a 83, 186 75 Praha 8, Czech Republic}\\
{\small E-mail: kalenda@karlin.mff.cuni.cz}
\and
{\sc Wies{\l}aw Kubi\'s}
\footnote{Research supported by MNiSW grant No. N201 024 32/0904.}
\\
{\small Institute of Mathematics of the Academy of Sciences of the Czech Republic}\\
{\small \v Zitn\'a 25, 115 67 Praha 1, Czech Republic}\\
\textit{and}\\
{\small Department of Mathematics, Jan Kochanowski University}\\
{\small \'Swi\c etokrzyska 15, 25-406 Kielce, Poland}
}

\begin{document}
\maketitle

\begin{abstract}
We study linearly ordered spaces which are Valdivia compact in their order topology. We find an internal characterization of these spaces and we present a counter-example 
disproving a conjecture posed earlier by the first author.
The conjecture asserted that a compact line is Valdivia compact if its weight does not exceed $\aleph_1$, every point of uncountable character is isolated from one side and every closed first countable subspace is metrizable. It turns out that the last condition is not sufficient. On the other hand, we show that the conjecture is valid if the closure of the set of points of uncountable character is scattered. This improves an earlier result of the first author.

\ 




\end{abstract}


\section{Introduction}

By a {\em compact line} we mean a linearly ordered compact space, i.e. a compact space whose topology is induced by a linear order.
We investigate Valdivia compact lines, i.e., compact lines which are Valdivia compact spaces. Recall that a compact space $K$ is called {\it Valdivia} if
it is homeomorphic to some $K'\subs\Err^\Gamma$ for a set $\Gamma$ such that
$$\{x\in K' : \{\gamma\in\Gamma: x(\gamma)\ne\emptyset \}\mbox{ is countable}\}$$
is dense in $K'$. Valdivia compact spaces play an important role in the study of the structure of nonseparable Banach spaces. They appeared first in \cite{AMN},
the name was given in \cite{DG}. For a detailed study of this class we refer to
\cite{Kalenda2000, valexa}.
Valdivia compact lines were addressed in \cite[Section 5]{K2006a}, \cite{K2006} and in \cite[Section 3]{valexa}. In \cite{valexa} the following question was asked.

\begin{question}\label{q-conj} Let $K$ be a compact line satisfying the following three conditions.
\begin{itemize}
	\item[(i)] $K$ has weight at most $\aleph_1$.
	\item[(ii)] Each point $x\in K$ of uncountable character is isolated from one side (i.e, one of the intervals $(\leftarrow,x]$ or $[x,\to)$ is open in $K$).
	\item[(iii)] Each closed first countable subset of $K$ is metrizable. 
\end{itemize}
Is $K$ necessarily Valdivia?
\end{question}

It is showed there that these conditions are necessary and that they are also sufficient in case $K$ is either scattered or connected. It is conjectured there \cite[Conjecture 3.5]{valexa}
that these conditions are sufficient in general.  In the present paper we show that the conjecture is false (Example~\ref{sanfowi} below). We further give a characterization of zero-dimensional Valdivia compact lines using functions defined on stationary subsets of $\omega_1$ (Theorem~\ref{0dimchar}). In this characterization we use a strengthening of condition (iii) from the above question. A characterization of general, not necessarily zero-dimensional, Valdivia compact lines is given in Section~\ref{sec-gencase}. In Section~\ref{sec-spec} we show that the above question has positive answer if the points of uncountable character have scattered
closure in $K$, which generalizes the results of \cite{valexa}. In the last section we study compact lines which are continuous images of Valdivia compacta.

\section{Preliminaries}

In this section we collect some auxiliary results on Valdivia compacta and namely on Valdivia compact lines, needed in the sequel.

A {\em compact line} is a compact space $K$ together with a linear order which induces the topology of $K$. It is well known that a linearly ordered set $X$ is compact in its interval topology if and only if it is order complete, i.e. every nonempty subset of $X$ has both the supremum and the infimum. A compact line $K$ is zero-dimensional if and only if for every $x,y\in K$ with $x<y$ there are $x',y'$ such that $x\loe x'<y'\loe y$ and the open interval $(x',y')$ is empty.
Given a compact line $K$, we shall denote by $0_K$ and $1_K$ the minimal and the maximal element of $K$ respectively.

We shall use standard notation concerning intervals in a linearly ordered set. For example: $[a,\rightarrow)$ will denote the closed final interval (segment) induced by $a$, i.e. $[a,\rightarrow)=\setof{x}{a\loe x}$.
A subset $G$ of a linearly ordered set $X$ is {\em convex} if $[x,y]\subs G$ whenever $x,y\in G$ are such that $x<y$. The smallest convex set containing $A\subs X$ will be denoted by $\conv (A)$.
A map $\map fXY$ between linearly ordered sets is {\em increasing} if $f(x_0)\loe f(x_1)$ whenever $x_0\loe x_1$. We shall often use the simple fact that every increasing surjection between compact lines is continuous.

We treat ordinals as well ordered sets with respect to $\in$. In particular, given two ordinals $\al,\beta$, $\al<\beta$ holds iff $\al\in\beta$.
Recall that $\omega_1$ denotes the smallest uncountable ordinal, which is at the same time treated as a linearly ordered space, endowed with the order topology. We shall denote by $\omega_1^{-1}$ the set $\omega_1$ with reversed ordering. Note that a set $C\subs \omega_1$ is {\em unbounded} if it has cardinality $\aleph_1$. Recall that a set $S\subs\omega_1$ is {\em stationary} if it intersects every closed unbounded subset of $\omega_1$. The {\em Pressing Down Lemma} says that given a stationary set $S\subs\omega_1$, for every function $\map fS{\omega_1}$ which is {\em regressive}, i.e. $f(\al)<\al$ for $\al\in S$, there exists a stationary set $S'\subs S$ on which $f$ is constant. For more information concerning ordinals and set-theoretic notions we refer to \cite{Jech} and \cite{Kunen}.

An {\em $\omega_1$-sequence} in a topological space $X$ is a function $\map x{\omega_1}X$. We shall often write $x_\al$ instead of $x(\al)$.
The notion of a limit of an $\omega_1$-sequence $x$ is defined naturally. Namely, $p=\lim_{\al\to\omega_1}x_\al$ if for every neighborhood $U$ of $p$ there is $\al<\omega_1$ such that $\setof{x_\xi}{\xi\goe\al}\subs U$. If $X$ is linearly ordered and the sequence $x$ is increasing then its possible limit is the supremum of the set $\setof{x_\al}{\al<\omega_1}$. A sequence is {\em monotone} if it is either increasing or decreasing (i.e. increasing with respect to the reversed ordering). 

Let $K$ be a compact space and $A\subs K$. We say that $A$ is a {\em $\Sigma$-subset} of $K$ if there is a homeomorphic injection $h:K\to\Err^\Gamma$ such that $A=h^{-1}(\Sigma(\Gamma))$,
where
$$\Sigma(\Gamma)=\setof{x\in\Err^\Gamma}{\{\gamma\in\Gamma: x(\gamma)\ne\emptyset \}\mbox{ is countable}}.$$
Hence, $K$ is Valdivia if and only if it admits a dense $\Sigma$-subspace.
	
Further, if $K$ is a compact line, following \cite{valexa} we denote by $G(K)$ the set of all points of $K$ which are either isolated or can be obtained as the limit of a one-to-one sequence. The following lemma was proved in \cite[Lemma 3.1]{valexa}.

\begin{lm}\label{G(K)} 
Let $K$ be a compact line. Then $G(K)$ is dense in $K$. Moreover, if $K$ is Valdivia, then
$G(K)$ is the unique dense $\Sigma$-subset of $K$ and is formed by all $G_\delta$ points of $K$.
\end{lm}

We will also use the following characterization of Valdivia compact lines.

\begin{lm}\label{separfamily}
Let $K$ be a compact line. Then $K$ is Valdivia if and only if there is a family $\Aaa$ of open $F_\sigma$-intervals in $K$ satisfying the following two conditions:
\begin{itemize}
	\item Family $\Aaa$ separates points of $K$, i.e. for each distinct points $x,y\in K$ there is $I\in\Aaa$ containing exactly one of them.
	\item Each $x\in G(K)$ belongs to countably many elements of $\Aaa$.
\end{itemize}
Further, if $\Aaa$ is such a family, we have
$$G(K)=\{x\in K:\{I\in\Aaa:x\in I\}\mbox{ is countable}\}.$$
If $K$ is moreover zero-dimensional, the family $\Aaa$ may be chosen to consist of clopen intervals.
\end{lm}

\begin{pf} It follows from Lemma~\ref{subsetsoflines} and \cite[Proposition 1.9]{Kalenda2000} that being Valdivia is equivalent to the existence of a family of open $F_\sigma$ sets satisfying the two conditions. We have also the equality given above. Finally, as each open $F_\sigma$ subset of a compact space is Lindel\"of, it can be expressed as a countable union of elements of a given basis. Therefore $\Aaa$ may be chosen to consist of open $F_\sigma$-intervals and in case $K$ is zero-dimensional to consist of clopen intervals.
\end{pf}

The next result is proved in \cite[Proposition 3.2]{valexa} and shows that Valdivia compact lines have rather exceptional structure.

\begin{lm}\label{subsetsoflines}
Let $K$ be a compact line. If $K$ is Valdivia (a continuous image of a Valdivia compact space),
then so is each closed subset of $K$.
\end{lm}

We will further need the following result on continuous images:

\begin{lm}\label{imagesoflines}
Let $K$ be a Valdivia compact line, $L$ a compact line and $\varphi:K\to L$ an order-preserving continuous surjection. Suppose that for each $y\in L$ either $\varphi^{-1}(y)$ is a singleton or $\varphi^{-1}(y)\cap G(K)$ is dense in $\varphi^{-1}(y)$. Then $L$ is Valdivia as well.
\end{lm}

\begin{pf} Set
$$E=\{\pair{x_1}{x_2}\in K\times K: \varphi(x_1)=\varphi(x_2)\}.$$
As $G(K)$ is a dense $\Sigma$-subset of $K$, by \cite[Lemma 2.9 and Theorem 2.20]{valchar}
it is enough to observe that $E\cap (G(K)\times G(K))$ is dense in $E$. But this easily follows from the assumptions. 
\end{pf}

\section{Zero-dimensional Valdivia compact lines}

In this section we give an internal characterization of zero-dimensional Valdivia compact lines.
We first restrict to the $0$-dimensional case as there is a duality between compact $0$-dimensional lines and linearly ordered sets. Namely, given a $0$-dimensional compact line $K$, let $\xk K$ be the set of all clopen final segments $F$ of $K$ such that $0_K\notin F$ and $1_K\in F$. Then $\xk K$ is a linearly ordered set (the order being defined by inverse inclusion). Conversely, given a linearly ordered set $X$, let $\kx X$ be the set of all final segments endowed with the topology inherited from the Cantor cube $\{0,1\}^X$, where each final segment is identified with its characteristic function. Then $\kx X$ is a compact $0$-dimensional line, the order is given by inverse inclusion. 
Note that $\kx{\xk{K}}$ is canonically order-homeomorphic to $K$ for each zero-dimensional compact line $K$ and $\xk{\kx{X}}$ is canonically order-isomorphic to $X$ for each linearly ordered set $X$.

The above defined operations naturally extend to contravariant functors which witness the isomorphism between the category of linearly ordered sets with increasing maps and the category of nonempty compact $0$-dimensional lines $K$ with continuous increasing maps.

The promised characterization of zero-dimensional Valdivia compact lines is
contained in the following theorem.

\begin{tw}\label{0dimchar} Let $X$ be a linearly ordered set. Then $\kx X$ is Valdivia if and only if the following three conditions are satisfied.
\begin{enumerate}
	\item[(1)] $|X|\loe\aleph_1$.
	\item[(2)] Every bounded monotone $\omega_1$-sequence has a limit in $X$.
	\item[(3)] For every stationary set $S\subs\omega_1$ and every map $f:S\to X$ there is a stationary set $T\subs S$ such that $f|_T$ is monotone.
\end{enumerate}
\end{tw}

Let us comment a bit the conditions in the above theorem. 

As the cardinality of $X$ is equal to the weight of $\kx X$, condition (1) means just that the weight of $\kx X$ is at most equal to $\aleph_1$. This corresponds to condition (i) in Question~\ref{q-conj}.

Condition (2) formulated in more detail means that each increasing $\omega_1$-sequence which is bounded from above has a supremum in $X$ and each decreasing $\omega_1$-sequence bounded from below has an infimum in $X$. Supposing that (1) holds, the validity of (2) is equivalent to the validity of condition (ii) from Question~\ref{q-conj} for the space $\kx X$. Indeed, if, say, $(x_\alpha)_{\alpha<\omega_1}$ is an increasing $\omega_1$-sequence which is bounded from above but has no supremum, then the final segment
$$\bigcap_{\alpha<\omega_1}(x_\alpha,\to)$$
has uncountable character in $\kx X$ while it is not isolated from either side.
Conversely, suppose that $k\in \kx X$ has uncountable character and is not isolated from either side. Without loss of generality we can suppose that the character of $k$ in $(\leftarrow,k]$ is uncountable. Then there are $k_\alpha\in (\leftarrow,k)$, $k<\omega_1$, isolated from the left such that the $\omega_1$-sequence $(k_\alpha)$ is increasing and has limit $k$. If we set
$x_\alpha=[k_\alpha,\to)$, we get an incresing $\omega_1$-sequence in $X$ which is bounded from above and having no limit in $X$.

As we will see below, condition (3) is a natural strengthening of condition (iii) from Question~\ref{q-conj}. 
We first prove the necessity of a weaker assumption.

\begin{prop}\label{pik_qow}
Let $X$ be a linearly ordered set. If\/ $\kx X$ is Valdivia compact then the following condition is satisfied:
\begin{enumerate}
	\item[(3')] Every uncountable subset of $X$ contains either a copy of $\omega_1$ or a copy of $\omega_1^{-1}$.
	\end{enumerate}
\end{prop}

\begin{pf}
Assume $Y\subs X$ contains neither $\omega_1$ nor $\omega_1^{-1}$. Then $\kx Y$ is a first countable increasing quotient of $\kx X$, therefore it is Corson compact by the result of \cite{Kalenda2000}. Nakhmanson's theorem \cite{Na} implies that $\kx Y$ is metrizable, therefore $|Y|\loe\aleph_0$.
\end{pf}

Note that condition (3') is weaker than (3). Indeed, suppose that (3) holds for a linearly ordered set $X$. Let $Y\subs X$ be uncountable. Then there is a one-to-one map $f:\omega_1\to Y$. By (3) there is a stationary set $T\subs\omega_1$ such that $f|_T$ is monotone.  If $f|_T$ is increasing,
then $f[T]$ is a copy of $\omega_1$, otherwise it is a copy of $\omega_1^{-1}$.

Further, note that condition (3') implies the validity of (iii) for the space $\kx X$. Indeed, let $L\subs \kx X$ be a closed first countable set. As $\kx X$ is zero-dimensional, there is an increasing retraction $r:\kx X\to L$. It follows that $\xk L$ is order-isomorphic to a subset of $X$. As $L$ is first countable, $\xk L$ contains no copies of $\omega_1$ or $\omega_1^{-1}$. By (3') we get that $\xk L$ is countable, so $L$ is metrizable.

Below (in Example~\ref{sanfowi}) we show that conditions (1), (2) and (3') are not sufficient for $\kx X$ being Valdivia. In particular, this will disprove the above conjecture. But before proving the example we need two lemmata.

\begin{lm}\label{jpwjefe}
Let $X, Y$ be linearly ordered sets such that $Y\subs X$. Denote by $f$ the dual map to the inclusion (hence $f$ is an incresing quotient mapping of $\kx X$ onto $\kx Y$). The following conditions are equivalent:
\begin{enumerate}
	\item[(a)] The  map $f$ is topologically right-invertible.
	\item[(b)] There exists an increasing map $\map p{\conv(Y)}Y$ such that $p|_Y = \id Y$.
	\item[(c)] Every proper gap $\pair AB$ in $Y$ remains a gap in $X$, i.e., 
	whenever $A,B\subs Y$ are nonempty subsets such that $A\cup B=Y$, $a<b$ whenever $a\in A$ and $b\in B$, $A$ has no maximum  and $B$ has no minimum, then there is no $x\in X$ with $a<x<b$ for all $a\in A$ and $b\in B$.
	\end{enumerate}
\end{lm}

\begin{pf}
(a)$\Rightarrow$(b) Let $g:\kx Y\to\kx X$ be a right inverse of $f$, i.e., $f\cmp g=\id {\kx Y}$. It is clear that $g$ must be increasing. We define $p$ as follows. Let $y\in\conv(Y)$. Set $k=(y,\to)$. Then $k\in \kx X$ and $k>0_{\kx X}$. As $y\in\conv(Y)$, there are $y_1,y_2\in Y$ with $y_1\loe y\loe y_2$.
Then
$$0_{\kx Y}\le [y_1,\to)\cap Y<(y_1,\to)\cap Y\loe (y,\to)\cap Y=f(k),$$
hence $k>g(0_{\kx Y})$. Further, given any $k_1<k$ we have
$$f(k_1)\loe [y,\to)\cap Y\loe[y_2\to)\cap Y<(y_2,\to)\cap Y\le 1_{\kx Y},$$
hence $k\loe g(1_{\kx Y})$.
Thus $[k,\to)\cap g[\kx Y]$ is a clopen final segment of $g[\kx Y]$ not containing $g(0_{\kx Y})$ but containing $g(1_{\kx Y})$. So, we can set
$p(y)=f[[k,\to)]$. It is now clear that $p$ maps $\conv(Y)$ into $Y$ and that it is increasing. Finally if $y\in Y$, then $[y,\to)\cap Y<(y,\to)\cap Y$ and hence $f[[k,\to)\cap Y]=[f(k),\to)$. The latter clopen interval corresponds to $y$. 
This completes the argument.

(b)$\Rightarrow$(c) Let $p$ be such a mapping. Let $A,B\subs Y$ be like in (c).
Suppose there is $x\in X$ such that $a<x<b$ for each $a\in A$ and $b\in B$.
Then $x\in\conv(Y)$ and so $p(x)$ is defined. As $Y=A\cup B$ necessarily $p(x)\in A$ or $p(x)\in B$. If $p(x)\in A$, then $a=p(a)\loe p(x)$ for each $a\in A$, so $p(x)$ is the maximum of $A$, a contradiction. Similarly, if $p(x)\in B$, then $p(x)$ is the minimum of $B$, a contradiction.

(c)$\Rightarrow$(a) We will define a right inverse of $f$ as follows. Take 
$k\in\kx Y$. Then $f^{-1}(k)$ is a closed interval in $\kx X$, say $[\alpha,\beta]$. If it is a singleton, the definition of $g(k)$ is clear.
Otherwise, if $\beta$ is isolated from the right, set $g(k)=\alpha$. If $\beta$ is not isolated from the right but $\alpha$ is isolated from the left, set $g(k)=\beta$. If this can be done for each $k\in\kx Y$, it is clear that $g$ is the required right inverse.

It remains to show that it is not possible that $\alpha<\beta$, $\alpha$ is not isolated from the left and $\beta$ is not isolated from the right. Suppose this possibility takes place. Let $A$ be the set of all elements of $Y$ such that the respective clopen interval has minimum less than $k$ and $B$ be the set of all elements of $Y$ such that the respective clopen interval has minimum greater than $k$. Then $A$ and $B$ do satisfy all assumptions given in (c). As $\alpha<\beta$, and $\kx X$ is zero-dimensional, there is a clopen interval in $\kx X$ with minimum in $(\alpha,\beta]$. The corresponding element of $x$ produces a contradiction with (c). 
\end{pf}

\begin{lm}\label{PRI}
Let $X$ be a linearly ordered space. Then $\kx X$ is Valdivia if and only if
there is a family $(X_\alpha:\alpha<\omega_1)$ satisfying the following properties:
\begin{itemize}
	\item[(i)] $X_\alpha$ is a countable subset of $X$;
	\item[(ii)] $X_\alpha\subs X_\beta$ for $\alpha<\beta$;
	\item[(iii)] $X_{\lambda}=\bigcup_{\alpha<\lambda} X_\alpha$ for every limit ordinal $\lambda<\omega_1$;
	\item[(iv)] $X=\bigcup_{\alpha<\omega_1} X_\alpha$;
	\item[(v)] the inclusion $X_\alpha\subs X$ satisfies condition (c) of Lemma~\ref{jpwjefe}.
\end{itemize}
\end{lm}

\begin{pf} It follows from \cite[Proposition 2.6 and Corollary 4.3]{kub-mich} that a compact space $K$ of weight $\aleph_1$ is Valdivia if and only if there is an $\omega_1$-sequence of retractions $(r_\alpha:\alpha<\omega_1)$
satisfying
\begin{itemize}
	\item $r_\alpha[K]$ is metrizable for each $\alpha<\omega_1$;
	\item $r_\alpha\cmp r_\beta=r_\beta\cmp r_\alpha=r_\alpha$ for $\alpha\le\beta<\omega_1$;
	\item the map $\alpha\mapsto r_\alpha(x)$ is continuous (when $\omega_1$ is equipped with the order topology) and has limit $x$  for each $x\in K$.
\end{itemize}

Moreover, if $K$ is linearly ordered, the retractions may be chosen increasing (this follows using \cite[Proposition 5.7]{K2006}).

If $K=\kx X$ is Valdivia, take such retractions and set $X_\alpha=\xk {r_\alpha[K]}$ canonically embedded into $X$. Then the family $(X_\alpha)$ satisfies the required conditions, the last one follows from Lemma~\ref{jpwjefe}.

Conversely, let $(X_\alpha)$ be a family satisfying the required conditions.
Set $X_{\omega_1}=X$ and $K_\alpha=\kx{X_\alpha}$ for $\alpha\le\omega_1$.
If $\alpha\le\beta\le\omega_1$ let $f^\beta_\alpha:K_\beta\to K_\alpha$ be the incresing surjection dual to the inclusion $X_\alpha\subs X_\beta$. Then it is clear that $f^\beta_\alpha=f^\gamma_\alpha\cmp f^\beta_\gamma$ for $\alpha\le\gamma\le\beta\le\omega_1$, hence we have an inverse sequence indexed by $\omega_1$. Moreover, this sequence is continuous because of condition (iii), all bonding maps are right-invertible by Lemma~\ref{jpwjefe}, $K_\alpha$ is metrizable for $\alpha<\omega_1$, therefore the limit $K=K_{\omega_1}$ is Valdivia compact by \cite[Corollary 4.3]{kub-mich}.
\end{pf}

Below is the announced example. In fact, it is a classical construction due to Kurepa \cite{Kurepa}, generalized by Todor\v cevi\'c in \cite[Section 4]{To}.

\begin{ex}\label{sanfowi} There is a linearly ordered set $Z$ satisfying conditions (1), (2) and (3') such that $\kx Z$ is not Valdivia.
\end{ex}

\begin{pf}
Let $X=\setof{x\in\Qyu^{\omega_1}}{|\suppt(x)|<\aleph_0}$ be endowed with the lexicographic ordering, where $\suppt(x) = \setof{\al}{x(\al)\ne0}$. Then $\kx X$ is a Valdivia compact. Indeed, the sets
$$X_\alpha=\setof{x\in X}{\suppt x\subset\alpha},\qquad \alpha<\omega_1,$$
have all properties from Lemma~\ref{PRI}.

We shall now extend $X$ by adding some new elements.
Fix a set $S\subs\omega_1$ consisting of limit ordinals. For each $\delta\in S$ choose a set $c_\delta$ order isomorphic to $\omega$ and such that $\sup(c_\delta)=\delta$. Now let
$$Y_S=\setof{1_{c_\delta}}{\delta\in S},$$
where $1_a$ denotes the characteristic function of the set $a\subs\omega_1$.
Define $X_S=X\cup Y_S$.
Clearly, $X_S$ satisfies (1). We check that $X_S$ satisfies (2).

We will use the following easy observation: If $\{a_\alpha\}_{\alpha<\omega_1}$ is a monotone $\omega_1$-sequence in $\Qyu^{\omega_1}$, then
\begin{equation*}(\forall\; \gamma<\omega_1)(\exists\;\alpha_0<\omega_1) \; \; \{a_\alpha\rest \gamma\}_{\alpha_0\le\alpha<\omega_1} \mbox{ is constant}.
\tag{*}\label{wdefowef}
\end{equation*}

Fix a strictly monotone sequence $\setof{a_\al}{\al<\omega_1}\subs X_S$. Define
$$T=\setof{t\in\Qyu^{<\omega_1}}{\exists\;\al<\omega_1\;\forall\;\xi\goe\al\;\;(t\subs a_\xi)}.$$
By $\Qyu^{<\omega_1}$ we mean $\bigcup_{\alpha<\omega_1}\Qyu^\alpha$, i.e. functions with rationals values whose domain is a countable ordinal. We consider $\Qyu^{\omega_1}$ ordered by inclusion. 
Note that $T$ is a chain in $\Qyu^{<\omega_1}$. Let $g=\bigcup T$. Then either $g\in T$ or $\dom(g)=\omega_1$. The first possibility cannot occur, because assuming $\delta=\dom(g)<\omega_1$ we would find (due to (\ref{wdefowef})) $\al_0<\omega_1$ such that $a_\al \rest\delta+1$ is constant for $\al\ge\al_0$ and then $a_{\al_0}\rest \delta+1 = g\cup\sn{\pair\delta {a_{\al_0}(\delta)}}$ would be an element of $T$.

Thus $g\in \Qyu^{\omega_1}$. It is clear that $\suppt(g)$ finite, because the sequence is strictly monotone and hence it contains at most one of the added elements $1_{c_\delta}$, $\delta\in S$. Thus $g\in X\subs X_S$. Further, we will show that $g$ is the limit of $\sett{a_\al}{\al<\omega_1}$ in $\Qyu^{\omega_1}$. 

Suppose that the sequence $\sett{a_\al}{\al<\omega_1}$ is increasing. Then $a_\al\le g$ for all $\alpha<\omega_1$. Indeed, suppose that there is some $\alpha_0<\omega_1$ with $a_{\alpha_0}>g$. Then for each $\alpha\ge\alpha_0$ we have $a_\alpha>g$ and so there is some $\gamma(\alpha)<\omega_1$ such that
$a_\alpha\rest \gamma(\alpha)=g\rest \gamma(\alpha)$ and $a_\alpha(\gamma(\alpha))>g(\gamma(\alpha))$. As $\sett{a_\al}{\al<\omega_1}$ is increasing, the $\omega_1$-sequence  $\sett{\gamma(\al)}{\alpha_0\ge \al<\omega_1}$ is decreasing. Therefore it is eventually constant, i.e., there is $\alpha_1\in[\alpha_0,\omega_1)$ and $\gamma<\omega_1$ such that for each $\alpha\in[\alpha_1,\omega_1)$ we have $\gamma(\alpha)=\gamma$. It follows that
$g\rest(\gamma+1)\notin T$, a contradiction.

Finally, it follows easily from the definition of $g$ that it is the supremum of 
$\sett{a_\al}{\al<\omega_1}$. If $\sett{a_\al}{\al<\omega_1}$ is decreasing, the proof is similar.
 This completes the proof of (2).

We now show that $X_S$ satisfies (3').
We shall use the fact that $\kx X$ is Valdivia compact. In particular, by Proposition \ref{pik_qow}, every uncountable subset of $X$ contains a copy of $\omega_1$ or $\omega_1^{-1}$. 
Fix (if possible) an uncountable set $A\subs S$. Let us denote $y_\delta=1_{c_\delta}$ for $\delta\in S$. We shall show that $\setof{y_\delta}{\delta\in A}$ contains a monotone subsequence.
For each limit ordinal $\lambda<\omega_1$ fix $\delta(\lambda)\in A$ such that $\delta(\lambda)>\lambda$ and set $\gamma(\lambda)=\sup\suppt (y_{\delta(\lambda)}\rest \lambda)$. By the Pressing Down Lemma there is a stationary set $S$ and $\gamma_0<\omega_1$ such that $\gamma(\lambda)=\gamma_0$ for each $\lambda\in S$. Moreover, as $\gamma_0$ is countable and $\suppt (y_{\delta(\lambda)}\rest\lambda)$ is finite for each $\lambda$, there are a stationary set $S'\subs S$ and a finite set $F\subs \gamma_0+1$
such that $\suppt (y_{\delta(\lambda)}\rest\lambda)=F$ for each $\lambda\in S'$.

Further we choose $\lambda_\eta\in S'$ for $\eta<\omega_1$ such that $\lambda_{\eta}>\delta(\lambda_{\theta})$ whenever $\theta<\eta<\omega_1$.
It can be done as $S'$ is unbounded in $\omega_1$. Finally note that $\sett{y_{\delta(\lambda_\eta)}}{\eta<\omega_1}$ is decreasing by the definition of the lexicographic order. This finishes the proof of (3').

Finally, notice that $\kx{X_S}$ is Valdivia compact if and only if $S$ is not stationary. This follows from Lemma~\ref{PRI}. Indeed, set 
$$X_\delta=\{x\in X_S: \exists \gamma<\delta: \suppt x\subs\gamma\}.$$
Then the family $(X_\delta:\delta<\omega_1)$ satisfies all conditions from Lemma~\ref{PRI} except for the last one. If there is a closed unbounded set $C\subs\omega_1\setminus S$, then the family $(X_\delta:\delta\in C)$ witnesses
that $\kx{X_S}$ is Valdivia. 

Conversely, suppose that $S$ is stationary and $\kx{X_S}$ is Valdivia. Let $(Y_\delta:\delta<\omega_1)$ be the family witnessing it (i.e., satisfying all the conditions from Lemma~\ref{PRI}). Now, there is a closed unbounded set $C\subs \omega_1$ such that $X_\delta=Y_\delta$ for each $\delta\in C$. Choose some $\delta\in C\cap S$. Then the sets
\begin{align*}
A&=\{x\in X_\delta:x<y_\delta\}\\
B&=\{x\in X_\delta: x>y_\delta\}
\end{align*}
witness that condition (c) of Lemma~\ref{jpwjefe} is violated for $Y_\delta\subs X_S$. This is a contradiction showing that $\kx{X_S}$ is not Valdivia.
\end{pf}

Before proving Theorem~\ref{0dimchar}, we briefly recall the method of elementary substructures which we use here.

In what follows, the letter $\chi$ will denote an uncountable regular cardinal, big enough so that all relevant objects have cardinality strictly less than $\chi$. More precisely, denote by $H(\chi)$ the family of all sets $x$ whose transitive closure $\tc x$ has cardinality $<\chi$. Recall that $\tc x = x\cup\bigcup x\cup\bigcup\bigcup x \cup \dots$. Now, saying ``$\chi$ is big enough" means that all objects under consideration (e.g. a given topological space, a given transformation, etc.) belong to $H(\chi)$.

The structure $\pair{H(\chi)}\in$ satisfies all the axioms of set theory, except possibly the power set axiom. A set $M$ is an {\em elementary substructure} of of $\pair{H(\chi)}\in$ if $M\subs H(\chi)$ and for every formula $\phi(x_1,\dots,x_n)$, for every $a_1,\dots,a_n\in M$, $M\models \phi(a_1,\dots,a_n)$ if and only if $H(\chi)\models\phi(a_1,\dots,a_n)$. Here, ``$M\models \phi$" means ``$M$ satisfies $\phi$" in the usual sense of model theory. The method of elementary submodels is based on the well known L\"owenheim-Skolem Theorem, saying that every countable subset of $H(\chi)$ can be enlarged to a countable elementary substructure of $H(\chi)$. 
As a consequence, given a countable set $S\subs H(\chi)$, one can easily consrtuct by induction a chain $\sett{M_\al}{\al<\omega_1}$ of countable elementary substructures of $\pair{H(\chi)}\in$ such that $S\subs M_0$ and $\al\subs M_\al$ for every $\al<\omega_1$. In fact, we may even require that $M_\al\in M_{\al+1}$ and that the chain be continuous, i.e. $M_\delta=\bigcup_{\xi<\delta}M_\xi$ for every limit ordinal $\delta$. The last property follows from the fact that $\bigcup_{\xi<\delta}M_\xi$ is again elementary. Given such a chain $\sett{M_\al}{\al<\omega_1}$ and setting $\delta_\al=M_\al\cap\omega_1$, we note that each $\delta_\al$ is a countable ordinal and the set $C=\setof{\delta_\al}{\al<\omega_1}$ is closed and unbounded in $\omega_1$. Consequently, if $S$ is a stationary subset of $\omega_1$, there exists $\al$ such that $M_\al\cap \omega_1\in S$. We shall use this remark below.

We refer to \cite{kub-mich, K2006a} for applications of elementary submodels in the context of retractions and Valdivia compacta. More explanations of the method and its use for finding projections in Banach spaces can be found in \cite{K_ps}.
Last but not least, \cite{Dow} is an important survey on the use of elementary substructures in general topology.

\begin{pf}[Proof of Theorem~\ref{0dimchar}]
Suppose first that $\kx X$ is Valdivia. Then (1) and (2) are satisfied by the above remarks. Let us prove (3).
Let $(X_\alpha:\alpha<\omega_1)$ be a family with properties from Lemma~\ref{PRI}. Let $\map {p_\al}{\conv(X_\al)}{X_\al}$ be an increasing projection, i.e. $p_\al\rest X_\al=\id{X_\al}$ (it exists by Lemma~\ref{jpwjefe}). Let us extend $p_\al$ by setting $p_\al(x)=-\infty$ if $x<\conv(X_\al)$ and $p_\al(x)=+\infty$ if $x>\conv(X_\al)$. Assuming $-\infty<x<+\infty$ for every $x\in X$, this defines an increasing map from $X$ into $X_\al\cup\dn{-\infty}{+\infty}$. We shall write $y_\al$ instead of $f(\al)$.

Fix a sufficiently big regular cardinal $\chi$ and fix a continuous chain $\sett{M_\al}{\al<\omega_1}$ of elementary substructures of $\pair{H(\chi)}{\in}$ such that $X\in M_0$, $f\in M_0$,  $\sett{p_\al}{\al<\omega_1}\in M_0$ and $\al\subs M_\al$ for $\alpha<\omega_1$.

If $\beta<\alpha<\omega_1$, then $\beta\in M_\alpha$, so $p_\beta\in M_\alpha$ and hence $X_\beta\in M_\alpha$ (as $X_\beta$ is the range of $p_\beta$). As $X_\beta$ is countable, we get $X_\beta\subset M_\alpha$ by \cite[Proposition 2]{K_ps}. Therefore $X\subset \bigcup_{\alpha<\omega_1} M_\alpha$, and so
$$C_1=\{\alpha<\omega_1: X_\alpha=X\cap M_\alpha\}$$
is a closed unbounded set.

Let $\delta_\al=\omega_1\cap M_\al$ and let $C_2=\setof{\al<\omega_1}{\delta_\al=\al}$. Then $C_2$ is a closed unbounded subset of $\omega_1$, so $C_1\cap C_2\cap S$ is stationary. 
Note that each $\delta_\al$ is a limit ordinal, therefore $p_{\al}(y_\al)\in X_{\xi(\al)}\cup\dn{-\infty}{+\infty}$ for some $\xi(\al)<\al$.
Using the Pressing Down Lemma, we may assume that $\xi(\al)=\xi$ for $\al\in S'$, where $S'\subs C_1\cap C_2\cap S$ is stationary.

Now suppose that for a stationary set $S_1\subs S'$ we have that $p_\al(y_\al)=-\infty$. Then the sequence $\setof{y_\al}{\al\in S_1}$ is strictly decreasing. Indeed, let $\alpha,\beta\in S_1$ such that $\alpha<\beta$.
Then $\alpha\in M_\beta$, so $y_\alpha\in M_\beta$. As $\alpha\in C_1$, we get
$y_\alpha\in X_\beta$. Further, $p_\beta(y_\beta)=-\infty$ and so $y_\beta<\conv(X_\beta)$. In particular $y_\beta<y_\alpha$.

Similarly, if the set $S_2=\setof{\al\in S'}{p_\al(y_\al)=+\infty}$ is stationary, we get a strictly increasing sequence $y\rest {S_2}$.
So assume that the set $$S''=\setof{\al\in S'}{p_\al(y_\al)\in X_\xi}$$ is stationary.
Using the fact that $X_\xi$ is countable, further refining $S''$ we may assume that $p_\al(y_\al) = v\in X_\xi$ for all $\al \in S''$. 
Now suppose $\al,\beta\in S''$ are such that $\al<\beta$ and $v<y_{\al}$ and $v<y_\beta$.
Then $p_\beta(y_\beta)=v$ and $p_\beta(y_\al)=y_\al$, because $y_\al\in X_\beta$. Since $p_\beta$ is order preserving, necessarily $y_\beta<y_\al$. This observation shows that $y\rest R$ is strictly decreasing, where
$$R=\setof{\al\in S''}{v<y_\al}.$$
Similarly, $y\rest L$ is strictly increasing, where
$$L=\setof{\al\in S''}{y_\al<v}.$$
One of these sets must be stationary, unless $y$ has constant value $v$ on a stationary set. This completes the proof of (3).

Now we are going to prove sufficiency. Let $X$ satisfy conditions (1)--(3).
Write $X=\bigcup_{\al<\omega_1}X_\al$, where $\wek x=\sett{X_\al}{\al<\omega_1}$ is an increasing chain of countable subsets of $X$ such that $X_\delta=\bigcup_{\xi<\delta}X_\xi$ whenever $\delta$ is a limit ordinal.

Denote by $S$ the set of all ordinals $\al<\omega_1$ for which there exist a proper gap $\pair{A_\al}{B_\al}$ in $X_\al$ and an element $y_\al\in X\setminus X_\al$ which fills this gap, i.e. $a<y_\al<b$ whenever $a\in A_\al$, $b\in B_\al$.
If there exists a closed unbounded set $C\subs \omega_1\setminus S$ then we are done by Lemma~\ref{PRI}. So suppose $S$ is stationary. Using (3), we fix a stationary subset $T\subs S$ such that $\wek y=\sett{y_\al}{\al\in T}$ is monotone. Note that $\wek y$ cannot be constant, because $y_\al\notin X_\al$ for $\al\in T$. Thus, reversing the order if necessary, we may assume that $\wek y$ is strictly increasing.

Fix a big enough regular cardinal $\chi$ and  fix a continuous chain $\sett{M_\al}{\al<\omega_1}$ of elementary substructures of $\pair{H(\chi)}{\in}$ such that $\wek x\in M_0$, $\wek y\in M_0$ and $\alpha\subset M_\alpha$ for $\alpha<\omega_1$. Similarly as above we get that
$X_\alpha\subset M_\beta$ for $\alpha<\beta<\omega_1$ and hence 
$$C_1=\{\alpha<\omega_1: X_\alpha=X\cap M_\alpha\}$$
is a closed unbounded set. 

Denote again $\delta_\al=\omega_1\cap M_\al$ and let $C_2=\setof{\al<\omega_1}{\delta_\al=\al}$. Then $C_2$ is a closed unbounded subset of $\omega_1$ and so is $C_1\cap C_2$.

Fix $\delta\in C_1\cap C_2\cap T$.
 Clearly, $\delta$ is a limit ordinal, so $X_\delta=\bigcup_{\xi<\delta}X_\xi$. Recalling that $y_\delta$ fills the gap $\pair{A_\delta}{B_\delta}$ we see that
\begin{equation}
M_\delta \models (\exists\;b\in X)(\forall\;\al\in T)\; y_\al<b.\tag{**}\label{werhoweqwr}
\end{equation}
To show it first note that $T\in M_0\subset M_\delta$, as $T$ is the domain of $\wek y$ and $\wek y\in M_0$. Further, we know that $X_\alpha\subset M_\delta$
for each $\alpha<\delta$, so $X_\delta\subset M_\delta$ as well. In particular,
$B_\delta\subset M_\delta$. So, choose some $b\in B_\delta$. Then $b\in M$.
Moreover, if $\al\in T\cap M$, then $\al<\delta$, so $y_\alpha<y_\delta<b$.
This proves (\ref{werhoweqwr}). By elementarity, the sequence $\wek y$ is bounded from above.

By (2) there exists $g\in X$ such that $g=\sup_{\al\in T}y_\al$. Find $\gamma\in C_1\cap C_2\cap T$ such that $\gamma\ge\delta$ and $g\in X_\gamma=M_\gamma\cap X$.

Observe that $[y_\gamma,g)\cap X_\gamma=\emptyset$. Indeed, if $x\in X_\gamma\cap [y_\gamma,g)$ then there would exist $\al\in T$ such that $x<y_\al<g$; by elementarity, we would have $x<y_\beta$ for some $\beta\in T\cap M_\gamma$ and hence $x<y_\gamma$, a contradiction.

Recalling that $g\in X_\gamma$, it follows that
$$g=\min \Bigl([y_\gamma,\rightarrow)\cap X_\gamma\Bigr)=\min B_\gamma.$$
This contradicts the fact that $y_\gamma$ fills the gap $\pair{A_\gamma}{B_\gamma}$.
\end{pf}

\section{The non-zero-dimensional case}\label{sec-gencase}

In this section we give a characterization of not necessarily zero-dimensional
Valdivia compact lines. 

Let $K$ be a Valdivia compact line. We introduce an equivalence relation $\sim$ on $K$ by setting $x\sim y$ if the interval $[x,y]$ is connected. It is clear that this is indeed an equivalence relation and that equivalence classes are closed intervals. As a closed subset of a Valdivia compact line is again Valdivia by Lemma~\ref{subsetsoflines}, each equivalence class is a connected Valdivia compact line. By \cite[Theorem 5.2]{K2006a} there are only five such spaces.
Two metrizable ones -- the singleton and the unit interval $[0,1]$ and three others, which are denoted by $R^\to+1$, $(R^{\to}+1)^{-1}$ and $^{\leftarrow}I^\to$ in \cite{K2006a}. $R^{\to}$ denotes the long line, i.e. the lexicographic product $\omega_1\cdot [0,1)$, $R^\to+1$ is its compactification made by adding the endpoint. The space $(R^\to+1)^{-1}$ is the order inverse of $R^\to+1$. Finally, the space $^\leftarrow I^\to$ is the unique connected compact line $[a,b]$ such that $a<b$ and for each $y\in(a,b)$ the interval $[y,b]$ is order homeomorphic to $R^\to+1$ and $[a,y]$ is order homeomorphic to $(R^{\to}+1)^{-1}$. Therefore each of the equivalence classes of $\sim$ is order isomorphic to one of these five spaces.

Further, as $K$ has weight at most equal to $\aleph_1$, at most $\aleph_1$ equivalence classes contain more than one point. Finally, the space
$$K_0=K\setminus\bigcup\{(a,b): a\sim b\}$$
is Valdivia as well (by Lemma~\ref{subsetsoflines}). Moreover, it is zero-dimensional, hence the criterion from the previous section applies. The space $K$ must also satisfy conditions (i)--(iii) from Question~\ref{q-conj}. 
Hence we have proved the necessity in the following theorem.

\begin{tw}\label{non0dim} Let $K$ be a compact line. Define the equivalence $\sim$ as above and define $K_0$ by the above formula. The space $K$ is Valdivia if and only if the following conditions are satisfied.
\begin{itemize}
	\item Each equivalence class is a Valdivia compact.
  \item The space $K_0$ is Valdivia.
  \item Each point of uncountable character is isolated from one side.
\end{itemize}
\end{tw}

\begin{pf} It remains to prove the sufficiency. We will use Lemma~\ref{separfamily}. Suppose the above three conditions are satisfied. Let $\Aaa$ be a family of clopen intervals in $K_0$ such that $\Aaa$ separates points of $K_0$ and for each $x\in G(K_0)$ there are only countably many elements of $\Aaa$ containing $x$.

For each interval $I\in\Aaa$ we define an open $F_\sigma$ interval $\widetilde I\subs K$ as follows. As $I$ is clopen, we have $I=[a,b]$ where $a$ is isolated from the left and $b$ is isolated from the right (in $K_0$).
If there is some $x<a$ such that $x\sim a$, choose some $a'\in(x,a)$.
Otherwise set $a'=a$. Similarly, if there is some $y>b$ with $y\sim b$, choose $b'\in(b,y)$. Otherwise set $b'=b$. Now set
$$\widetilde I=\begin{cases} [a,b],& \mbox{ if }a'=a,b'=b,\\
                             (a',b],& \mbox{ if }a'<a,b'=b,\\
                             [a,b'),& \mbox{ if }a'=a,b'>b,\\
                             (a',b'),& \mbox{ if }a'<a,b'>b.\end{cases}$$
It is clear that $\widetilde I$ is an open $F_\sigma$ interval in $K$. Indeed, if $a'=a$, then $a$ is isolated from the left also in $K$, and if $a'<a$, then  $a'$ has countable character in $K$ (as it is not isolated from either side).
Similarly for $b$ and $b'$. 

Set $\widetilde \Aaa=\{\widetilde I: I\in \Aaa\}$. Then $\widetilde \Aaa$ separates points of $K_0$ and for each $x\in G(K)$ there are only countably many elements of $\widetilde\Aaa$ containing $x$.

Indeed, if $x\in G(K)\cap K_0$, then $x\in G(K_0)$ and hence it is only in countably many elements of $\Aaa$. As $\widetilde I\cap K_0=I$ for each $I\in \Aaa$, $x$ belongs to only countably many elements of $\widetilde\Aaa$ as well.

Furher, suppose that $x\in G(K)\setminus K_0$. Let $[a,b]$ be the equivalence class containing $x$. Then necessarily $a,b\in G(K_0)$. (If, say, $a$ has uncountable character in $K_0$, then it is not isolated from either side in $K$, a contradiction.) Finally observe that if $x\in \tilde I$ for some $I\in \Aaa$, then necessarily either $a\in I$ or $b\in I$. Thus there can be at most countably many such $I$'s.

Now we extend $\widetilde \Aaa$ in order to separate points of $K$. Fix any equivalence class $[a,b]$ with $a<b$. Then $a,b\in K_0$. There are four possibilities:

(i) $[a,b]$ is order-homeomorphic to $[0,1]$. Then let $\Bee_{a,b}$ be a countable basis of $(a,b)$ consisting of intervals.

(ii) $[a,b]$ is order-homeomorphic to $^\leftarrow I^\to$. Then $[a,b]$ is clopen in $K$. Let $\Bee_{a,b}$ be a family of open $F_\sigma$ intervals in $[a,b]$ separating points of $[a,b]$ such that $[a,b]\in \Bee_{a,b}$ and each point of $(a,b)$ is contained only in countably many elements of $\Bee_{a,b}$.

(iii) $[a,b]$ is order-homeomorphic to $R^\to+1$. Then $b$ is isolated from the right in $K$. Let $\Cee$  be a family of open $F_\sigma$ intervals in $[a,b]$ separating points of $[a,b]$ such that $[a,b]\in \Cee$ and each point of $[a,b)$ is contained only in countably many elements of $\Cee$. Set $\Bee_{a,b}=\{I\cap(a,b]:I\in\Cee\}$.

(iv) $[a,b]$ is order-homeomorphic to $(R^\to+1)^{-1}$. Define the family $\Bee_{a,b}$ analogously as in (iii).

Finally set
$$\Yu=\widetilde\Aaa\cup\bigcup\{\Bee_{a,b}: a,b\in K_0,a<b,a\sim b\}.$$
Then $\Yu$ is a family of open $F_\sigma$ intervals in $K$ separating points of $K$ such that each point of $G(K)$ is contained only in countably many elements of $\Yu$. This proves that $K$ is Valdivia.
\end{pf}

\section{When the conjecture is valid}\label{sec-spec}

In this section we prove the following theorem.

\begin{tw} Let $K$ be a compact line satisfying the conditions (i)--(iii) from Question~\ref{q-conj} such that the closure of the set of points of uncountable character is scattered. Then $K$ is Valdivia.
\end{tw}

\begin{pf} We define an equivalence relation on $K$ by setting $a\sim b$ if and only if $[a,b]$ is Valdivia. We will prove, using the conditions (i)--(iii), that it is really an equivalence relation and that the quotient $L=K/\sim$ is a connected compact line. Further, if we denote by $q$ the canonical quotient mapping of $K$ onto $L$, we will show that the image under $q$ of the set of points of uncountable character is dense in $L$ unless $L$ is a singleton. 
If we prove all this, it follows that $L$ is a connected compact line which is at the same time scattered (as a continuous image of a scattered compact space is again scattered). Thus $L$ is a singleton, hence $K$ is Valdivia.

So, it is enough to prove the above mentioned properties of $\sim$.

\paragraph{Step 1.} $\sim$ is an equivalence relation.

First note that by Lemma~\ref{subsetsoflines} we have $x\sim y$ whenever $x,y\in[a,b]$ and $a\sim b$. To prove that it is really an equivalence relation it remains to check that $a\sim c$ whenever $a<b<c$ and $a\sim b$ and $b\sim c$.
If $b$ is isolated from one side, then $[a,c]$ is the topological sum of two Valdivia compact spaces and therefore it is Valdivia as well. If $b$ is isolated from neither side, it has countable characted. As the topological sum $[a,b]\oplus[b,c]$ is Valdivia and $[a,c]$ is the quotient made by identifying the two copies of $b$, we get that $[a,c]$ is Valdivia by Lemma~\ref{imagesoflines}.

\paragraph{Step 2.} Each equivalence class of $\sim$ is closed. 

Let $a\in K$ be arbitrary and let $b$ be the supremum of all $x\in K$ with $a\sim x$. We will show that $a\sim b$. Suppose $b>a$. 

If $b$ is isolated from the left, then clearly $a\sim b$ (in this case the supremum in the definition of $b$ is obviously attained).

So suppose that $b$ is not isolated from the left. If there is some $x\in(a,b)$ such that $[x,b]$ is connected, then $[x,b]$ is a connected compact line satisfying (i)--(iii) and hence it is Valdivia by \cite[Section 3.2]{valexa}, so $x\sim b$. As $a\sim x$, we get $a\sim b$.

Next suppose that $b$ is not isolated from the left and $[x,b]$ is connected for no $x\in (a,b)$. There are two possibilities:

(a) There is a sequence $b_n$ in $(a,b)$ such that $b_n\nearrow b$. Without loss of generality we may suppose that each $b_n$ is isolated from one side. By passing to a subsequence we may suppose that all the points $b_n$ are isolated from the same side. Suppose that each $b_n$ is isolated from the left. Then $[a,b]$ is homeomorphic to the one-point compactification of the topological sum of Valdivia compacta $[a,b_1)$, $[b_1,b_2), [b_2,b_3),\dots$.
Therefore $[a,b]$ is Valdivia by \cite[Theorem 3.35]{Kalenda2000} and so $a\sim b$. If each $b_n$ is isolated from the right we can proceed similarly.

(b) $b$ has uncountable character in $[a,b]$. Then we can find an increasing homeomorphic copy $(b_\alpha:\alpha<\omega_1)$ of $\omega_1$ in $(a,b)$ with supremum $b$. Again, we can without loss of generality suppose that for each isolated ordinal $\alpha<\omega_1$ the point $b_\alpha$ is isolated from one side and that all these points are isolated from the same side. Suppose they are isolated from the left. Set 
$X_\alpha=[b_{\alpha-1},b_{\alpha})$ for $\alpha<\omega_1$ isolated (where $b_{-1}=a$. Let $X$ be the $[0,\omega_1+1)$-sum of these spaces in the sense of \cite{weakval}. Then $X$ is Valdivia by \cite[Proposition 3.4]{weakval}. Moreover, $X$ is homeomorphic to the compact line made from $[a,b]$ by duplicating $b_\alpha$ for each limit ordinal $\alpha<\omega_1$. As each of these points has countable character, by Lemma~\ref{imagesoflines} we get that $[a,b]$ is Valdivia and so $a\sim b$.
If  $b_\alpha$ is isolated from the right for each isolated $\alpha<\omega_1$, we can proceed similarly.

\paragraph{Step 3.} $L=K/\sim$ is a connected compact line.

By Steps 1 and 2 we get that the quotient $L=K/\sim$ is a Hausdorff compact space, in fact a compact line. Let us prove that $L$ is connected. Suppose that $a,b\in L$ are such that $a<b$ and $(a,b)=\emptyset$. 
Let $[x_1,y_1]$ be the equivalence class in $K$ corresponding to $a$ and $[x_2,y_2]$ be the equivalence class corresponding to $b$. Then $y_1<x_2$ and $[y_1,x_2]=\{y_1,x_2\}$ and therefore $y_1\sim x_2$, a contradiction.

\paragraph{Step 4.} Denote by $q$ the quotient mapping of $K$ onto $L$ associated to $\sim$. Then the image under $q$ of the set of all points of uncountable character in $K$ is dense in $L$ unless $L$ is a singleton. 

To see this, let $y\in L$ be any point different from the endpoint. The inverse image of $y$ is a closed interval $[a,b]\subs K$. If $b$ is isolated from the right, then there is $b^+=\min(b,\to)$. Clearly $b\sim b^+$, hence also $a\sim b^+$. This is a contradiction to the fact that $[a,b]$ is an equivalence class of $\sim$. Further, let $c>b$ be arbitrary. If $[b,c]$ is first countable, then it is metrizable by (iii) and so $b\sim c$. Again we get $a\sim c$ which is a contradiction. Thus $[b,c]$ contains a point of uncoutable character. Therefore there is a net $d_\nu$ of points of uncountable character converging to $b$. Finally, $q(d_\nu)\to y$.

The proof is complete.
\end{pf}

\section{Compact lines which are continuous images of Valdivia compacta}

In this section we discuss the class of compact lines which are continuous images of Valdivia compacta. The main question in this context is the following one.

\begin{question}\label{q-image} Let $K$ be a compact line which is a continuous image of a Valdivia compact space. Is there a Valdivia compact line $L$ and an order-preserving continuous surjection of $L$ onto $K$?
\end{question}

In \cite{valexa} it is conjectured that the answer is positive (Conjecture 3.6).
It follows from \cite[Theorem 3.7]{valexa} that the answer is positive for scattered compact lines. Let us note that the weight of a compact line which is a continuous image of a Valdivia compact cannot exceed $\aleph_1$. This fact can be proved by similar arguments as in \cite[Prop. 5.5]{K2006a}.

We do not know the answer to this question. Nonetheless, we provide a characterization of order-preserving quotients of Valdivia compact lines.

\begin{tw} Let $K$ be a compact line. Denote by $A$ the set of all points of uncountable character in $K$ which are not isolated from any side. Then $K$ is an order-preserving quotient of a Valdivia compact line if and only if the compact line made from $K$ by duplicating all points of $A$ is Valdivia.
\end{tw}

\begin{pf} The if part is obvious (the order preserving quotient is made by collating back the duplicated points). 

Let us show the only if part. Let $L$ be a Valdivia compact line and $\varphi:L\to K$ an order-preserving continuous surjection.
Let $\sim$ be the associated equivalence relation on $L$, i.e. $x\sim y$ if and only if $\varphi(x)=\varphi(y)$. By our assumptions the equivalence classes
are closed intervals. Set
$$L_1=L\setminus\bigcup\{(a,b):a\sim b\}.$$
Then $L_1$ is again a Valdivia compact line by Lemma~\ref{subsetsoflines}.
Moreover, $\varphi(L_1)=K$ and $\varphi$ is at most two-to-one. Denote the restriction of $\sim$ to $L_1$ again by $\sim$. Then each eqivalence class has at most two points.

For each equivalence class $\{a,b\}$ such that $a<b$ and at least one of the points $a,b$ is isolated in $L_1$, choose an isolated point $y_{a,b}\in\{a,b\}$.
Denote by $L_2$ the compact line made from $L_1$ by omitting all these points $y_{a,b}$. Then $L_2$ is again a Valdivia compact line (Lemma~\ref{subsetsoflines}) and $\varphi(L_2)=K$. Denote the restriction of $\sim$ to $L_2$ again by $\sim$. Then the equivalence classes have at most two points and, moreover, if an equivalence class has two points, none of them is isolated in $L_2$.

Define an equivalence relation $\sim_2$ on $L_2$ by the following formula
$$x\sim_2 y \Leftrightarrow x=y \mbox{ or } (x\sim y \mbox{ and both points $x,y$ have countable character in }L_2)$$
Then $L_3=L_2/\sim_2$ is a Valdivia compact line by Lemma~\ref{imagesoflines}.

Finally it is easy to see that $L_3$ is exactly the compact line made from $K$ by duplicating all points of uncountable character which are not isolated from any side.
\end{pf}

As a consequence we get that the non-Valdivia compact lines constructed in Example~\ref{sanfowi} are not order-preserving quotients of a Valdivia compact line.
Are they continuous images of a Valdivia compact space?

\subsection*{Acknowledgements}

The second author would like to thank the Department of Mathematical Analysis of Charles University in Prague for supporting his visits (Fall 2006, Fall 2007) when most of this work was being done.

\end{document}